\documentclass[12pt]{amsart}
\usepackage{amssymb,latexsym}
\usepackage{amsfonts}
\usepackage{amsmath}
\usepackage[colorlinks,linkcolor=blue,anchorcolor=blue,citecolor=blue]{hyperref}
\usepackage{algorithm}
\usepackage{enumerate}
\usepackage{algpseudocode}
\usepackage{verbatim}
\usepackage{graphicx}

\newcommand\C{{\mathbb C}}

\newcommand\Q{\mathbb{Q}}

\newcommand\Z{{\mathbb Z}}

\newcommand\oQ{\overline{\Q}}
\newcommand\N{\mathrm{(N)}}
\newcommand\PP{\mathrm{(P)}}

\newcommand\Gal{{\mathrm {Gal}}}

\newcommand\cA{\mathcal{A}}
\newcommand\cS{\mathcal{S}}

\newcommand\h{\mathrm{h}}
\newcommand\tor{\mathrm{tor}}
\newcommand\hh{\hat{\h}}
\newcommand\ab{\mathrm{ab}}
\newcommand\tr{\mathrm{tr}}
\newcommand\cc{\mathrm{c}}

\newcommand{\Gm}{\mathbb{G}_\mathrm{m}}
\newcommand\real{\mathrm{real}}

\newtheorem{theorem}{Theorem}[section]
\newtheorem{lemma}[theorem]{Lemma}

\newtheorem{corollary}[theorem]{Corollary}

\theoremstyle{definition}
\newtheorem{definition}[theorem]{Definition}

\newtheorem{question}[theorem]{Question}
\newtheorem{example}[theorem]{Example}
\newtheorem{conjecture}[theorem]{Conjecture}

\numberwithin{equation}{section}

\def\house#1{{%
    \setbox0=\hbox{$#1$}
    \vrule height \dimexpr\ht0+1.4pt width .5pt depth \dp0\relax
    \vrule height \dimexpr\ht0+1.4pt width \dimexpr\wd0+2pt depth \dimexpr-\ht0-1pt\relax
    \llap{$#1$\kern1pt}
    \vrule height \dimexpr\ht0+1.4pt width .5pt depth \dp0\relax}}

\begin{document}

\title[]{On the properties of Northcott and Narkiewicz for elliptic curves}

\author{Jorge Mello}
\address{Max Planck Institute for Mathematics, Bonn, 53111, Germany}
\email{jmelloguitar@gmail.com}

\author{Min Sha}
\address{School of Mathematical Sciences, South China Normal University, Guangzhou, 510631, China}
\email{shamin2021@qq.com}

\subjclass[2010]{11G05, 11G50}

\keywords{Elliptic curve, Northcott property, Property (P), height}

\begin{abstract}
In this paper, as an analogue of the number field case, for an elliptic curve $E$ defined over the algebraic numbers and for any subfield $F$ of algebraic numbers, 
we say that $E$ has the Northcott property over $F$ if there are at most finitely many $F$-rational points on $E$ of uniformly bounded height, 
and we say that $E$ has the property (P) over $F$ if for any infinite subset $S$ of $F$-rational points on $E$, 
$f(S) = S$ for an $F$-endomorphism $f$ of $E$ implies that $f$ is an automorphism. 
We establish some criteria for both properties and provide typical examples. 
We also show that the Northcott property implies the property (P), but the converse is not true.  
\end{abstract}

\maketitle

\section{Introduction}

\subsection{Background and motivation}

Let $\cA$ be a subset of the algebraic numbers $\oQ$. 
Following Bombieri and Zannier \cite{BZ}, we say that $\cA$ has the \textit{Northcott property}, short \textit{Property $\N$},  if for any positive number $T$ 
the set 
$$
\cA(T) = \{\alpha \in \cA: \, \h(\alpha) \le T \}
$$
is finite, where $\h(\alpha)$ denotes the absolute logarithmic Weil height (see \cite[Chapter VIII]{Silv}). 
The famous Northcott theorem \cite{North} states that a set of algebraic numbers of uniformly bounded degree 
(over the rational numbers $\Q$) has Property $\N$. 

Given a number field $K$ and a positive integer $d$, let $K^{(d)}$ be the compositum of all the field extensions of $K$ of degree at most $d$. 
Bombieri and Zannier \cite{BZ} asked whether the field $K^{(d)}$ has Property $\N$. 
So far, for this question, due to \cite[Corollary 1]{BZ} we only know that $K^{(2)}$ has Property $\N$. 
More generally, it was showed in \cite[Theorem 1]{BZ} that $K^{(d)}_{\ab}$ has Property $\N$ for any $d \ge 2$. 
Here, $K^{(d)}_{\ab}$ is the compositum of all the abelian extensions of $K$ contained in $K^{(d)}$. 
Another consequence is that for any positive integer $d$, the field $\Q(1^{1/d}, 2^{1/d}, 3^{1/d}, \ldots)$ 
has Property $\N$. 
Later, Dvornicich and Zannier \cite[Theorem 2.1]{DZ} proved that 
for any field $F$ of algebraic numbers with Property $\N$, 
any finite extension of $F$ also has Property $\N$. 
Recently, Widmer \cite[Theorem 3]{Wid} obtained a criterion for Property $\N$  
which yields new examples of fields having Property $\N$.  See \cite{CF,CW,Fe} for further work. 

A subfield $F$ of $\oQ$ is said to have the \textit{property} (P) if 
for any infinite subset $S$ of $F$, 
$f(S) = S$ for a polynomial $f$ over $F$ implies $\deg f = 1$. 
This property was first introduced and investigated by Narkiewicz in \cite{Nark}. 
Both \cite{DZ} and \cite{Nark} gave some typical examples for the property (P). 
Dvornicich and Zannier also showed in \cite[Theorem 3.1]{DZ} that the Northcott property implies the property (P). 
Pottmeyer proved in \cite[Corollary 4.5]{Pot1} that $\Q^{(d)}$ has the property (P) for any positive integer $d$. 
Recently, Pottmeyer \cite{Pot3} studied the property (P) in the field generated by all the symmetric Galois extensions of $\Q$. 
In addition, in \cite[Section 6]{CW} Checcoli and Widmer investigated the relations between the property (P) 
and other properties of fields; see also \cite[Section 4]{Pot1}. 

In this paper, we want to define and study Property $\N$ of subsets of 
 algebraic points on an elliptic curve defined over a number field, as well as relate it to a similar version of the property (P) in this context.

\subsection{The Northcott property}

Throughout the paper, $E$ is an elliptic curve defined over $\oQ$ by a Weierstrass equation, 
and $O$ is the point at infinity of $E$. 
For any subfield $F$ of $\oQ$, $E(F)$ denotes the set of $F$-rational points on $E$. 
As usual, when $E$ is defined over $F$,  $E(F)$ can be made to be an abelian group, and let $E(F)_\tor$ be the torsion subgroup of $E(F)$.  

\begin{definition}
For any subset $\cS$ of $E(\oQ)$, we say that $\cS$ has the \textit{Northcott property}, short \textit{Property $\N$}, 
if for any positive number $T$ the set 
$$
\cS(T) = \{P \in \cS: \, \h(P) \le T \}
$$
is finite, where the height $\h(P)$ is defined to be the height of its $x$-coordinate (see Section~\ref{sec:height}). 
If $E(F)$ has Property $\N$, we also say that $E$ has \textit{Property $\N$ over $F$}. 
\end{definition}

Certainly, the above notion can be extended to any algebraic variety defined over $\oQ$ 
after choosing a height function. 

We first show that Property $\N$ is preserved under non-zero isogeny. 

\begin{theorem}  \label{thm:isogeny}
Let $E_1, E_2$ be two elliptic curves  defined over $\oQ$, and let $F$ be a subfield of $\oQ$. 
 Assume that  $E_1$ and $E_2$ are isogenous over $F$. 
Then,  $E_1(F)$ has Property $\N$ if and only if $E_2(F)$ does.
\end{theorem}

The following is a consequence of Theorem~\ref{thm:isogeny} by replacing non-zero isogeny with non-constant morphism 
but imposing that the curves are also defined over $F$. 

\begin{corollary}  \label{cor:morphism}
Let $F$ be a subfield of $\oQ$, and let $E_1, E_2$ be two elliptic curves  defined over $F$. 
Assume that there is a non-constant morphism defined over $F$ from $E_1$ to $E_2$. 
Then, $E_1(F)$ has Property $\N$ 
if and only if $E_2(F)$ does.
\end{corollary}

If one could give a positive answer to the following question, 
then in Corollary~\ref{cor:morphism} we would not need to assume that the curves are both defined over $F$ 
 (because we could make deductions over a finite extension of $F$). 

\begin{question} \label{quest:ext}
Let $L/F$ be a finite extension of subfields of $\oQ$. 
If $E(F)$ has Property $\N$, does $E(L)$ also have Property $\N$? 
\end{question}

We try to answer Question~\ref{quest:ext} under some extra conditions. 

If $E$ is defined over $F$, then for any finite Galois extension $L/F$ of subfields of $\oQ$, we can define the trace map: 
$$
T_{L/F}: \, E(L) \to E(F), \quad  P \mapsto \sum_{\sigma \in \Gal(L/F)} P^\sigma, 
$$
where $P^\sigma = (\sigma(x),\sigma(y))$ if $P=(x,y)$. 

\begin{theorem} \label{thm:ext2}
Let $L/F$ be a finite Galois extension of subfields of $\oQ$. 
Assume that $E$ is defined over $F$ and $E(F)$ has Property $\N$. 
Then, $E(L)$ has Property $\N$ if and only if the kernel of the trace map $T_{L/F}$ does. 
\end{theorem}

Moreover, we provide a sufficient condition to get a positive answer for Question~\ref{quest:ext}.

\begin{theorem} \label{thm:ext4}
Let $E, F, L$ be as in Theorem~\ref{thm:ext2}. 
Assume further that the quotient group $E(L)/E(F)$ is finitely generated.  
Then, $E(L)$ also has Property $\N$. 
\end{theorem}

We remark that in general for a finite extension $L/F$ of subfields of $\oQ$, 
the quotient group $E(L)/E(F)$ is not always finitely generated. 
For example, let $E$ be defined by $y^2 = x^3 + x$ over $\Q$, and let $L = \Q^\cc$ 
the cyclotomic closure of $\Q$ and $F$ the maximal real subfield of $L$, 
then we have $[L:F] = 2$, but $E(L)/E(F)$ is not finitely generated due to the points $(- p, \sqrt{- p(p^2 + 1)})$ 
for primes $p$.

Now, we want to consider the following question. 

\begin{question} \label{quest:sub}
For which subfield $F$ of $\oQ$, $E(F)$ has Property $\N$?  
\end{question}

Note that the torsion points in $E(\oQ)$ are of  bounded height. 
So, a necessary condition for $E(F)$ having Property $\N$ is that 
$E(F)_\tor$ is finite. 

Clearly, if $F$ has Property $\N$, then so does $E(F)$. 
But the converse is not true in general (see Example~\ref{exam:Zp} below). 

Here we give a simple sufficient condition for a positive answer about Question~\ref{quest:sub}. 

\begin{theorem}  \label{thm:fg}
For any subfield $F$ of $\oQ$ and for any elliptic cruve $E$ defined over $F$, 
if the abelian group $E(F)$ is finitely generated, then $E(F)$ has Property $\N$. 
\end{theorem}

We remark that the converse of Theorem~\ref{thm:fg} is not true (see Example~\ref{exam:K2} below). 
However,  if the quotient group $E(F)/mE(F)$ is finite for some integer $m \ge 2$, 
then by the descent theorem (see \cite[Chapter VIII, Theorem 3.1]{Silv}), the converse of Theorem~\ref{thm:fg} is true. 

\begin{example}  \label{exam:K2}
Let $K$ be a number field and $F=K^{(2)}$. 
Then, for any elliptic curve $E$ defined over $K$, 
by the proof of \cite[Theorem 4.5]{Im} the rank of $E(F)$ is infinite. 
Since $F$ has Property $\N$ by Bombieri and Zannier, 
 $E(F)$ also has Property $\N$.
\end{example}

\begin{example}\label{exam:Zp}
Let $p_1, \ldots, p_n$ be distinct prime numbers. 
Denote $F = \Q(\mu_{p_1^{\infty}}, \ldots, \mu_{p_n^{\infty}})$, 
where, for any prime $p$, $\mu_{p^{\infty}}$ denotes the group of $p$-power roots of unity. 
Then, for any elliptic curve $E$ defined over $\Q$, 
$E(F)$ is finitely generated (this is a consequence of theorems of K. Kato, K. Ribet and D. Rohrlich; see \cite[Theorem 1.1]{Lozano-Robledo}), 
and so by Theorem~\ref{thm:fg}, $E(F)$ has Property $\N$. 
However, the field $F$ itself does not have Property $\N$, because $F$ contains infinitely many roots of unity.
\end{example}


Now, we want to establish some typical cases for the failure of Property $\N$. 

Let $K$ be a number field. That is, $K$ is a finite extension of $\Q$ and is a subfield of $\overline{\Q}$. Then, $\overline{\Q}$ is the algebraic closure of $K$. 
For any Galois isomorphisms $\sigma_1, \ldots, \sigma_n \in \Gal(\overline{\Q}/K)$, 
denote by $\overline{\Q}^{\sigma_1, \ldots, \sigma_n}$ the fixed subfield of $\overline{\Q}$ under $\sigma_1, \ldots, \sigma_n$. 
That is, 
$$
\overline{\Q}^{\sigma_1, \ldots, \sigma_n} = \{\alpha \in \overline{\Q}:\, \sigma_i(\alpha) = \alpha, i=1, \ldots, n\}. 
$$

\begin{theorem} \label{thm:Ksigman}
Let $E$ be an elliptic curve defined over a number field $K$. 
For any $\sigma_1, \ldots, \sigma_n \in \Gal(\overline{\Q}/K)$, 
if $E(\overline{\Q}^{\sigma_1, \ldots, \sigma_n})$ has infinite rank, then 
 $E(\overline{\Q}^{\sigma_1, \ldots, \sigma_n})$ does not have Property $\N$.
\end{theorem}

Im and Larsen showed in \cite[Theorem 1.5]{ImLar} that the set of $n$-tuples 
$(\sigma_1, \ldots, \sigma_n) \in \Gal(\overline{\Q}/K)^n$ for which
 $E(\overline{\Q}^{\sigma_1, \ldots, \sigma_n})$ has infinite rank 
 is of full measure. 
They also showed in \cite[Theorem 1.4]{ImLar} that 
for any $\sigma \in \Gal(\overline{\Q}/K)$,  $E(\overline{\Q}^\sigma)$ has infinite rank. 
Hence, we directly have:

\begin{corollary} \label{cor:Ksigma}
Let $E$ be an elliptic curve defined over a number field $K$. 
For any $\sigma \in \Gal(\overline{\Q}/K)$, 
 $E(\overline{\Q}^\sigma)$ does not have Property $\N$.
\end{corollary}

Let $\Q^\real$ be the field of the real algebraic numbers. 
The following corollary is a direct consequence of Corollary~\ref{cor:Ksigma} 
(choosing $\sigma$ to be the complex conjugation).

\begin{corollary} \label{cor:Qreal}
Let $E$ be an elliptic curve defined over $\Q^\real$. 
Then, $E(\Q^\real)$ does not have Property $\N$.
\end{corollary}

Let $\Q^\tr$ be the field of the totally real algebraic numbers. 
Note that $\Q^\tr$ is a proper subfield of $\Q^\real$. 

\begin{theorem} \label{thm:Qtr}
Let $E$ be an elliptic curve defined over $\Q^\tr$. 
Then, $E(\Q^\tr)$ does not have Property $\N$.
\end{theorem}

For any number field $K$, let $K^\cc$ be the cyclotomic closure of $K$, 
and let $K^\ab$ be the abelian closure of $K$. 
That is, $K^\cc$ is the field generated over $K$  by all the roots of unity, 
and $K^\ab$ is the compositum of all the abelian extensions of $K$. 
Clearly, $K^\cc$ is a subfield of $K^\ab$. 
By the Kronecker-Weber theorem, we have $\Q^\ab = \Q^\cc$. 
However, this is not true in general if $K \ne \Q$. 

Let $E$ be an elliptic curve defined over $K$. 
Ribet \cite{Ribet} showed that $E(K^\cc)_\tor$ is finite. 
Zarhin \cite{Zar} proved that $E(K^\ab)_\tor$ is finite if and only if $E$ does not have complex multiplication over $K$. 

Hence, if $E$ has complex multiplication over $K$, then $E(K^\ab)_\tor$ is infinite, 
and so $E(K^\ab)$ does not have Property $\N$. 
We conjecture that this is true in general. 

\begin{conjecture}  \label{conj:Kab}
Let $E$ be an elliptic curve defined over a number field $K$. 
Then,  $E(K^\ab)$ does not  have Property $\N$.  
\end{conjecture}

Finally, we present an application of Property $\N$, which suggests that Conjecture~\ref{conj:Kab} might be not easy. 

Let $E$ be an elliptic curve defined over $K^\cc$. 
By \cite[Corollary 1 (a)]{DZ07}, there are only finitely many points in $E(K^\cc)$ 
with a root of unity as $x$-coordinate or as $y$-coordinate. 
Here we generalise this result to elliptic curves defined over $\Q$ for  points whose coordinates are quotients of cyclotomic integers of bounded house 
(see Section~\ref{sec:cyclo} for the definition of house); see \cite{BS} for another generalisation. 
Here, cyclotomic integers mean algebraic integers in $\Q^\ab$. 

\begin{theorem} \label{thm:boundedhouse} 
Let $E$ be an elliptic curve defined by a Weierstrass equation over $\Q$. 
Then, there are only finitely many points in $E(\Q^\ab)$ whose coordinates are quotients of cyclotomic integers of uniformly bounded house.
\end{theorem}

\subsection{The property (P)} 
As an analogue of the field case \cite{Nark}, we define the property (P) for elliptic curves. 
We then relate it to the Northcott property as in \cite{DZ}. 

Recall that $E$ is an elliptic curve defined over $\oQ$. 

\begin{definition}
For any subset $\mathcal{S}$ of $E(\oQ)$, we say that $\mathcal{S}$ has the  \textit{property} (P)
if for any infinite subset $\Gamma \subset \mathcal{S}$, the condition $f(\Gamma)=\Gamma$ for  a endomorphism $f$ of $E$ defined over $\mathbb{Q}(\mathcal{S})$ implies $\deg f=1$ (that is, $f$ is an automorphism). 
Here $\mathbb{Q}(\mathcal{S})$ denotes the field generated by all the elements in $\mathcal{S}$ over $\Q$. 
\end{definition}

Differently from the Northcott property, a subset of $E(\oQ)$ having the property (P) might contain infinitely many torsion points.  
For example, if $E$ does not have complex multiplication, then, 
any subset $\cS$ of prime-order torsion points 
for which for any prime $p$, $\cS$ contains exactly one non-zero $p$-torsion point, 
 has the property (P). 

But in many cases finiteness of torsion points is indeed necessary. 

\begin{theorem}  \label{thm:Ptor}
Let $F$ be a subfield of $\overline{\Q}$. 
Assume that $E$ is defined over $F$ and $E(F)$ has the property $\PP$. 
Then, $E(F)$ contains only finitely many torsion points.  
\end{theorem}

As in \cite[Proposition 3.1]{DZ},  the property (P) implies finiteness of preperiodic points. 
Recall that given an endomorphism $f$ of $E$, a point $P$ on $E$ is called a periodic point of $f$ 
if $f^{(n)}(P) = P$ for some integer $n > 0$, where $f^{(n)}$ denotes the $n$-th iteration of $f$; 
and $P$ is called a preperiodic point of $f$ if $f^{(n)}(P)$ is a periodic point for some integer $n > 0$. 

\begin{theorem}\label{thm:preper}
Let $\cS$ be a subset of $E(\oQ)$, and let $f$ be an endomorphism of $E$ defined over $\Q(\cS)$ of degree greater than $1$. 
Assume that $\cS$ has the property $\PP$ and $f(\cS) \subseteq \cS$. 
Then, $f$ has only finitely many preperiodic points in $\cS$.
\end{theorem}

By Theorem~\ref{thm:preper}, we directly have: 
\begin{corollary}
If $E(F)$ has the property $\PP$ for a subfield $F \subseteq \oQ$, then every endomorphism $f$ of $E$ defined over $F$ of degree greater than $1$ 
has only finitely many preperiodic points in $E(F)$.
\end{corollary}

As the field case \cite[Theorem 3.1]{DZ}, we show that the Northcott property implies the property (P).

\begin{theorem}\label{thm:NimpliesP}
For any subset $\cS$ of $E(\oQ)$, the Northcott property for $\cS$ implies the property $\PP$ for $\cS$.
 \end{theorem}
 
From Theorem~\ref{thm:NimpliesP} we know that for any number field $K$, $E(K)$ has the property (P), 
because it has the Northcott property. 
 
 We remark that the converse of Theorem~\ref{thm:NimpliesP} is not true. 
 For this, we first establish a sufficient condition for the property (P). 
 
\begin{theorem}\label{thm:critP}
Let $\cS$ be a subset of $E(\oQ)$. 
Assume that $\cS$ has only finitely many torsion points and
the canonical heights of non-torsion points in $\cS$ have a uniform positive lower bound. 
Then, $\cS$ has the property $\PP$.
\end{theorem}

Combining Theorem~\ref{thm:critP} with \cite[Corollary 2]{Zhang}, we directly have:

\begin{corollary}
Let $E$ be an elliptic curve defined over a number field $K$.  Then, $E(K\Q^\tr)$ has the property $\PP$.
\end{corollary}

So, for any elliptic curve $E$ defined over $\Q^\tr$, $E(\Q^\tr)$ has the property (P), 
but by Theorem~\ref{thm:Qtr} $E(\Q^\tr)$ does not have the Northcott property. 
This gives a counterexample for the converse of Theorem~\ref{thm:NimpliesP}. 

Using Theorem~\ref{thm:critP}, we can obtain more examples for the property (P). 

Combining Theorem~\ref{thm:critP} with the result of 
Ribet \cite{Ribet}  mentioned before 
and together with \cite[Theorem 1]{Sil} we obtain: 

\begin{corollary} \label{cor:PK}
Let $E$ be an elliptic curve defined over a number field $K$. 
Then, $E(K^\cc)$ has the property $\PP$. 
\end{corollary}

Moreover, in Corollary~\ref{cor:PK} if $E$ does not have complex multiplication over $K$, 
using Zarhin's result \cite{Zar} instead of Ribet's we have that $E(K^\ab)$ has the property (P). 
However, if $E$ has complex multiplication over $K$,  
then by Zarhin's result \cite{Zar}, $E(K^\ab)$ contains infinitely many torsion points,  
and so by Theorem~\ref{thm:Ptor}, $E(K^\ab)$ does not have the property $\PP$. 

\begin{corollary} \label{cor:P-Kab}
Let $E$ be an elliptic curve defined over a number field $K$. 
Then, $E(K^\ab)$ has the property $\PP$ if $E$ does not have complex multiplication over $K$, 
and otherwise $E(K^\ab)$ does not have the property $\PP$.
\end{corollary}




Finally, \cite[Example 5.7]{Pot2} suggests that 
the property (P) is not always preserved under finite extensions (see below for more details).  

\begin{example}  \label{exam:pres}
We follow \cite[Example 5.7]{Pot2}. 
Let $E$ be an elliptic curve defined by $y^2 = x^3 + Ax + B$ over a number field $K$. 
Assume that $E$ has additive reduction at a finite place $v$ of $K$, and that 
there is an element $\gamma \in K$ such that the twist $E_{\gamma}$ 
(defined by $y^2 = x^3 + \gamma^2 Ax + \gamma^3 B$) has good reducton at $v$. 
Let $K^{nr, v}$ be the maximal algebraic extension of $K$ which is unramified above $v$. 
Then, by Theorem~\ref{thm:critP} and \cite[Proposition 5.4]{Pot2}, we know that 
$E(K^{nr, v})$ has the property (P). 
In addition, by \cite[Example 5.7]{Pot2}  (together with \cite[Proposition 5.6]{Pot2}), 
there is an infinite sequence of non-torsion points $Q_0, Q_1, \ldots$ in $E(K^{nr, v}(\sqrt{\gamma}))$ such that for some integer $m \ge 2$, 
$mQ_{n} = Q_{n-1}$ for any $n \ge 1$. 
Define $\Gamma = \{Q_0, Q_1, \ldots\} \cup \{mQ_0, m^2 Q_0, \ldots\}$. 
Clearly, $m\Gamma = \Gamma$. 
So, $E(K^{nr, v}(\sqrt{\gamma}))$ does not have the property (P). 
\end{example}

We will gather some preliminary results in Section~\ref{sec:pre} 
and then prove the main results one by one in Section~\ref{sec:proof}.

\section{Preliminaries}
\label{sec:pre}

\subsection{Height function} 
\label{sec:height}

For any non-zero algebraic number $\alpha$, $\h(\alpha)$ is the absolute logarithmic Weil height of $\alpha$. 
Let $E$ be an elliptic curve defined over $\overline{\Q}$. 
For any point $P \in E(\oQ)$, define the height of $P$ to be $\h(x(P))$, denoted by $\h_E(P)$, where $x(P)$ is the $x$-coordinate of $P$. 
The canonical height of $P$ is defined by 
$$
\hh_E(P) = \lim_{n \to \infty} \frac{\h_E(2^nP)}{4^n}. 
$$ 
For simplicity, we also denote $\h(P) = \h_E(P)$ and $\hh(P) = \hh_E(P)$ if the context is clear.  

One relation between $\h(P)$ and $\hh(P)$ is (see \cite[Chapter VIII, Theorem 9.3 (e)]{Silv})
\begin{equation*}   
\hh(P) = \h(P) + O(1),
\end{equation*}
where the $O(1)$ depends only on $E$. 
So, in defining Property $\N$, we can use $\hh(P)$ instead of $\h(P)$. 

Moreover,  for any non-zero endomorphism $f$ of $E$ of degree $d$, by \cite[Theorem B.2.5 (b)]{HS}, we have that 
for any $P \in E(\oQ)$,  
\begin{equation}   \label{eq:f-h}
 \h(f(P)) = d\h(P) + O(1), 
\end{equation}
where the implied constant depends only on $f$ and $E$. 
Then, noticing $2^n f(P) = f(2^n P)$ and by \eqref{eq:f-h}, we directly have 
\begin{equation}   \label{eq:f-hh}
 \hh(f(P)) = d\hh(P).
\end{equation}

In addition, for any $P\in E(\oQ)$ and $m\in \Z$, we have: 
\begin{itemize}
\item $\hh(mP)=m^2\hh(P)$; 
\item $\hh(P)=0$ if and only if $P$ is a torsion point.
\end{itemize}
Furthermore, for any $P,Q\in E(\oQ)$, we have 
\begin{equation} \label{eq:hhPQ}
\hh(P+Q) + \hh(P-Q) = 2\hh(P) + 2\hh(Q). 
\end{equation}

Following the hints in~\cite[Chapter~IX, Exercise~9.8]{Silv} and using~\cite[Chapter~VIII, Proposition~9.6]{Silv}, one can show that 
if $P_1,\ldots, P_r$ are linearly independent non-torsion points in $E(\oQ)$ (assuming $r\ge 1$), 
then for any integers $k_1,\ldots, k_r$, we have 
\begin{equation}  \label{eq:hhP1r}
\hh(k_1P_1 + \cdots + k_rP_r) \ge c\max_{1\le i \le r} k_i^2, 
\end{equation}
where $c$ is a positive constant depending only on $E$ and $P_1,\ldots, P_r$.

\subsection{Representations via sum of roots of unity} 
\label{sec:cyclo}

Let $\mu_\infty$ be the set of all the roots of unity in $\C$. 
For any algebraic number $\alpha$, 
 the \textit{house} of $\alpha$ is defined to be the maximum of the absolute values 
of all the conjugates of $\alpha$ over $\Q$.

Loxton~\cite[Theorem~1]{Loxton} showed that any cyclotomic integer $\alpha$ 
has a short representation as a sum of roots of unity, 
that is, $\alpha=\sum_{i=1}^m\zeta_i$, where $\zeta_1,\ldots,\zeta_m \in \mu_\infty$, and the integer $m$ depends only on
 the house of $\alpha$. 

We remark that Dvornicich and Zannier~\cite[Theorem~L]{DZ07} extended the above result to algebraic integers 
contained in a cyclotomic extension of a given number field. 



\subsection{Torsion points on varieties}

Let $n \ge 1$ be an integer. 
As usual, let $\Gm^n$ be the complex algebraic torus, that is, 
the $n$-fold Cartesian product of the multiplicative group $\Gm = \C^*$ of the complex numbers $\C$. 
The elements in $\mu_\infty^n$ are called the \textit{torsion points} of $\Gm^n$.

We say that an absolutely irreducible  variety $V \subseteq \Gm^n$ does not contain a \textit{monomial curve}, if it does not contain a curve parametrised by 
\begin{center}
 $x_1= \rho_1t^{k_1}, \, ..., \, x_n=\rho_n t^{k_n}$,
\end{center} 
where $\rho_1,..., \rho_n \in \mu_\infty$ and $k_1,...,k_n$ are integers, not all equal to zero.

The following is a special case of the toric Manin-Mumford conjecture proved by Laurent~\cite{Laurent}.

\begin{lemma}\label{lem:Laurent}
 If an algebraic variety $V$ in $\Gm^n$ does not contain a monomial curve, there are only finitely many torsion points of $\Gm^n$ on $V$.
\end{lemma}

\section{Proof of the main results}
\label{sec:proof}

\subsection{Proof of Theorem~\ref{thm:isogeny}}
Since $E_1$ and $E_2$ are isogenous over $F$, there is a non-zero isogeny defined over $F$, say $\phi$, from $E_1$ to $E_2$. 
Let $\hat{\phi}: E_2 \to E_1$ be the dual isogeny of $\phi$. 
Note that $\hat{\phi}$ is also defined over $F$. 

Since the $x$-coordinate function is an even function on $E_1$ and $\phi$ is in fact a group homomorphism, 
$x \circ \phi $ is also an even function on $E_1$, 
and therefore by~\cite[Chapter~VIII, Lemma~6.3]{Silv}, one has that for any $P \in E_1(\oQ)$, 
\begin{equation}\label{eq:h1h2}
(\deg \phi) \h_{E_1}(P) = \h_{E_2} \circ \phi(P) + O(1), 
\end{equation} 
where the implied constant does not depend on $P$. 

Suppose that $E_2(F)$ has Property $\N$.  
For any positive number $T$, define 
$$
S(T) = \{ P \in E_1(F) : \, \h_{E_1}(P) \leq T\}. 
$$
We need to show that $S(T)$ is a finite set. 
Note that we have (for instance, see \cite[Chapter II, Proposition 2.6 (a)]{Silv})
$$
 |S(T)|  \le (\deg \phi) |\phi(S(T))| .
$$
Clearly, $\phi(S(T))\subseteq E_2(F)$. 
By \eqref{eq:h1h2}, 
 the points in the set $\phi(S(T))$ are also of bounded height, 
 and so the set itself is a finite set by Property $\N$ of $E_2(F)$. 
Thus, $S(T)$ is also a finite set. 

Conversely, suppose that $E_1(F)$ has Property $\N$.  
For any positive number $T$, define 
$$
R(T) = \{ Q \in E_2(F)  : \, \h_{E_2}(Q) \leq T\}. 
$$
Note that we have 
$$
 |R(T)|  \le (\deg \hat{\phi}) |\hat{\phi}(R(T))| .
$$
As the above, we get that $R(T)$ is a finite set. 
This completes the proof.

\subsection{Proof of Corollary~\ref{cor:morphism}}
Let $\sigma: E_1 \to E_2$ be a non-constant morphism defined over $F$. 
Then, by \cite[Chapter~III, Example~4.7]{Silv}, there is a non-zero isogeny $\phi: E_1 \to E_2$  of the form  
$$
\phi = \tau \circ \sigma, 
$$
where $\tau$ is a translation map from $E_2$ to itself. 
Moreover, one can see that since both $E_1$ and $E_2$ are defined over $F$, $\tau$ is also defined over $F$. 
So, $\phi$ is in fact defined over $F$. 
That is, $E_1$ and $E_2$ are isogenous over $F$.
Now, the desired result follows from Theorem~\ref{thm:isogeny}.

\subsection{Proof of Theorem~\ref{thm:ext2}}
We only need to prove the sufficiency. 
By contradiction, we assume that $E(L)$ does not have Property $\N$. 
Then, there are infinitely many points, say $P_1, P_2, \ldots$, in $E(L)$ of uniformly bounded height. 
By \cite[Chapter VIII, Theorem 6.2]{Silv}, the points $T_{L/F}(P_i)$, $i=1,2, \ldots$, in $E(F)$ are also of uniformly bounded height. 
Since $E(F)$ has Property $\N$, there are at most finitely many distinct points in 
$T_{L/F}(P_1), T_{L/F}(P_2), \ldots$. 
For simplicity, we can assume that $T_{L/F}(P_1) = T_{L/F}(P_2) = \ldots$. 
Define $Q_i = P_i - P_1, i =2,3, \ldots$. 
Then, the points $Q_2, Q_3, \ldots$ have zero trace and are of uniformly bounded height. 
That is, they form an infinite set of points of uniformly bounded height in the kernel of $T_{L/F}$. 
This contradicts the hypothesis that the kernel of $T_{L/F}$ has Property $\N$.

\subsection{Proof of Theorem~\ref{thm:ext4}}
By Theorem~\ref{thm:ext2}, it suffices to show that the kernel of the trace map $T_{L/F}$, denoted by $\ker(T_{L/F})$,  has Property $\N$. 

Let $n = [L:F]$. 
By assumption, the quotient group $E(L)/E(F)$ is finitely generated. 
Then, the subgroup $\ker(T_{L/F})/E(F)$ is also finitely generated, say by the classes of $P_1, ..., P_r \in \ker(T_{L/F})$ modulo $E(F)$. 
So, for any $P \in \ker(T_{L/F})$, there are integers $k_1, ..., k_r$ such that 
$$
P \equiv k_1P_1+ \cdots + k_rP_r   \quad \mod{E(F)}, 
$$
and thus for some point $Q \in E(F)$ we have 
$$
P = k_1P_1+ \cdots + k_rP_r  + Q, 
$$
which, after taking the trace $T_{L/F}$, becomes  $nQ = O$, that is, $Q$ is an $n$-torsion point. 
Therefore, the kernel $\ker(T_{L/F})$ is finitely generated, 
and in particular it contains only finitely many torsion points. 

Without loss of generality, suppose that the points $P_1, P_2, \ldots, P_r$ form a basis 
of the free part of $\ker(T_{L/F})$. 
Then, every point $P \in \ker(T_{L/F})$ can be written as 
$$
P = k_1 P_1 + \cdots + k_r P_r + Q
$$
for some integers $k_1, \ldots, k_r$ and some torsion point $Q \in \ker(T_{L/F})$. 
Then, using \eqref{eq:hhPQ} and \eqref{eq:hhP1r} and noticing $\hh(Q)=0$,  we have 
\begin{align*}
\hh(P) & = \frac{1}{2} \left(2\hh(P) + 2\hh(Q)\right) = \frac{1}{2} \left(\hh(P+Q) + \hh(P-Q)\right) \\ 
& \ge \frac{1}{2} \hh(k_1 P_1 + \cdots + k_r P_r) \ge c \max_{1 \le i \le r} k_i^2, 
\end{align*}
where $c$ is a positive constant depending only on $E$ and $P_1, \ldots, P_r$. 
So, for any subset of points in $\ker(T_{L/F})$ of uniformly bounded height, 
the corresponding coefficients $k_1, \ldots, k_r$ are also uniformly bounded, 
and thus the subset itself is finite. 
This  completes the proof.

\subsection{Proof of Theorem~\ref{thm:fg}}
Applying the same arguements as the later part in the proof of Theorem~\ref{thm:ext4}, 
one can obtain the desired result.

\subsection{Proof of Theorem~\ref{thm:Ksigman}}
First, we claim that for any integer $m \ge 1$, we have 
\begin{equation} \label{eq:Ksigman}
\left| E(\overline{\Q}^{\sigma_1, \ldots, \sigma_n}) / m E(\overline{\Q}^{\sigma_1, \ldots, \sigma_n}) \right| \le m^{2n}. 
\end{equation}

We prove \eqref{eq:Ksigman} by induction on $n$. 
For $n=1$, we write $\sigma = \sigma_1$. 
By contradiction, suppose that $\left| E(\overline{\Q}^\sigma) / m E(\overline{\Q}^\sigma) \right| > m^2$. 
Then there exist points $P_1, \ldots, P_{k} \in E(\overline{\Q}^\sigma)$ with $k = m^2 +1$ 
such that for any $i \ne j$, $P_i - P_j \not\in m E(\overline{\Q}^\sigma)$. 
For each $1 \le i \le k$, choosing a point $Q_i \in E(\overline{\Q})$ with $P_i = mQ_i$, 
since $P_i^\sigma = P_i$, we have $m(Q_i^\sigma - Q_i) = O$, 
which means that $Q_i^\sigma - Q_i$ is an $m$-torsion point of $E$. 
Then, we obtain $k$ $m$-torsion points of $E$: $Q_1^\sigma - Q_1, \ldots, Q_k^\sigma - Q_k$. 
However, $E$ has exactly $m^2$ $m$-torsion points. 
Hence, there exist $1\le i \ne j \le k$ such that $Q_i^\sigma - Q_i = Q_j^\sigma - Q_j$, 
and so $(Q_i - Q_j)^\sigma = Q_i - Q_j$, which implies $Q_i - Q_j \in E(\overline{\Q}^\sigma)$. 
So, $P_i - P_j = m(Q_i - Q_j) \in mE(\overline{\Q}^\sigma)$, which contradicts with the choices of $P_1, \ldots, P_k$. 
Therefore, we obtain \eqref{eq:Ksigman} when $n=1$. 

For $n > 1$, by induction we assume that \eqref{eq:Ksigman} holds for $n-1$, that is, 
$$
\left| E(\overline{\Q}^{\sigma_1, \ldots, \sigma_{n-1}}) / m E(\overline{\Q}^{\sigma_1, \ldots, \sigma_{n-1}}) \right| \le m^{2(n-1)}.  
$$
Now, since $E(\overline{\Q}^{\sigma_1, \ldots, \sigma_n}) \subseteq E(\overline{\Q}^{\sigma_1, \ldots, \sigma_{n-1}})$, 
we consider the natural map from $E(\overline{\Q}^{\sigma_1, \ldots, \sigma_n}) / m E(\overline{\Q}^{\sigma_1, \ldots, \sigma_n})$ 
to $E(\overline{\Q}^{\sigma_1, \ldots, \sigma_{n-1}}) / m E(\overline{\Q}^{\sigma_1, \ldots, \sigma_{n-1}})$. 
Applying similar argument as in the above paragraph, we have that the kernel of this map 
has at most $m^2$ elements. 
Hence, we obtain \eqref{eq:Ksigman}. 

By contradiction, assume that $E(\overline{\Q}^{\sigma_1, \ldots, \sigma_n})$ has Property $\N$. 
Then, using \eqref{eq:Ksigman} and the descent theorem (see \cite[Chapter VIII, Theorem 3.1]{Silv}) 
we deduce that $E(\overline{\Q}^{\sigma_1, \ldots, \sigma_n})$ is finitely generated. 
However,  by assumption $E(\overline{\Q}^{\sigma_1, \ldots, \sigma_n})$ has infinite rank. 
We obtain a contradiction. 
Hence, $E(\overline{\Q}^{\sigma_1, \ldots, \sigma_n})$ does not have Property $\N$.

\subsection{Proof of Theorem~\ref{thm:Qtr}}
Using Corollary~\ref{cor:morphism}, we can assume that the elliptic curve $E$ is defined by a Weierstrass equation
$$
y^2 = x^3 + ax + b,  
$$
where $a, b \in \Q^\tr, 4a^3+27b^2 \ne 0$. 
For any integer $k \ge 1$, we define an elliptic curve $E_k$ by
$$
y^2 = (x+k)^3 + a(x+k) + b = x^3 + 3kx^2 + (3k^2+a)x + k^3+ak+b. 
$$
Note that there is an isomorphism from $E$ to $E_k$ defined by sending $(x,y)$ to $(x-k,y)$. 

Since both $a$ and $b$ have only finitely many Galois conjugates (they are all real numbers), 
we can choose a sufficiently large integer $k$ such that for any Galois isomorphism $\sigma$ we have 
\begin{equation}  \label{eq:kab}
k^3 + k \sigma(a) + \sigma(b) > 8 + 12k + 2(3k^2 + |\sigma(a)|). 
\end{equation}
Then, for any root of unity $\zeta$ and any Galois isomorphism $\sigma$, we have 
$$
(\zeta + \bar{\zeta})^3 + 3k(\zeta + \bar{\zeta})^2 + (3k^2+\sigma(a))(\zeta + \bar{\zeta}) + k^3+k\sigma(a)+\sigma(b) > 0,
$$
which in fact follows from \eqref{eq:kab} and 
$$
|(\zeta + \bar{\zeta})^3 + 3k(\zeta + \bar{\zeta})^2 + (3k^2+\sigma(a))(\zeta + \bar{\zeta})| \le 8 + 12k + 2(3k^2 + |\sigma(a)|). 
$$
Note that for any root of unity $\zeta$, $\zeta + \bar{\zeta}$ is a totally real number, where $\bar{\zeta}$ is the complex conjugate of $\zeta$. 
Hence, any point $(\zeta + \bar{\zeta}, \beta)$ on $E_k$ belongs to $E_k(\Q^\tr)$, 
and certainly all these points are of uniformly bounded height. 
Clearly, there are infinitely many distinct such points. 
So, $E_k(\Q^\tr)$ does not have Property $\N$. 
By Corollary~\ref{cor:morphism}, $E(\Q^\tr)$ also does not have Property $\N$.

\subsection{Proof of Theorem~\ref{thm:boundedhouse}}
Assume that $E$ is defined by the Weierstrass equation $f(x,y)=0$ over $\Q$. 

Without loss of generality (see the arguments below), we only need to prove that 
 there are only finitely many points in $E(\Q^\ab)$ whose coordinates are cyclotomic integers of uniformly bounded house.

By contradiction, suppose that there are infinitely many points on $E$ whose coordinates 
 are cyclotomic integers  with house bounded above by a certain number $B$. 
Then, by Loxton's result mentioned in Section~\ref{sec:cyclo}, there exists a positive integer $M=M(B)$ such that 
for any such point $(\alpha, \beta)$,  $\alpha$ and $\beta$ can be written as 
\begin{equation}\label{eq:albe}
\alpha= \sum_{i=1}^m  \zeta_i,\qquad \beta= \sum_{j=1}^n  \xi_j,
\end{equation} 
for some positive integers $m, n \le M$ and $\zeta_1, \ldots, \zeta_m, \xi_1, \ldots, \xi_n \in \mu_\infty$. 

Since by hypothesis there are infinitely many such points, there are fixed positive integers $m, n$ such that 
there are infinitely many points $(\alpha, \beta)$ on $E$ such that $\alpha$ and $\beta$ are cyclotomic integers 
and can be written as in \eqref{eq:albe}. 
This means that the algebraic variety $V$ defined by 
$$
f(x_1+ \cdots + x_m, y_1 + \cdots + y_n) = 0
$$ 
(with $x_1, \ldots, x_m, y_1, \ldots, y_n$ as variables) contains infinitely many torsion points of $\Gm^{m+n}$. 
By Lemma~\ref{lem:Laurent}, the variety $V$ contains  a monomial curve parametrised by 
$$
 x_1= \rho_1t^{k_1}, \ldots,  x_m=\rho_m t^{k_m},  y_1= \eta_1t^{l_1}, \ldots,  y_n=\eta_n t^{l_n}, 
$$
where $\rho_1, \ldots,  \rho_m, \eta_1, \ldots, \eta_n \in \mu_\infty$ and $k_1, \ldots, k_m, l_1, \ldots, l_n$ are integers, not all equal to zero. 
Then, we have 
$$
f(\rho_1t^{k_1}+ \cdots + \rho_m t^{k_m}, \eta_1t^{l_1} + \cdots + \eta_n t^{l_n}) = 0. 
$$ 
So, for any $\zeta \in \mu_\infty$, the point $(\rho_1\zeta^{k_1}+ \cdots + \rho_m \zeta^{k_m}, \eta_1\zeta^{l_1} + \cdots + \eta_n \zeta^{l_n})$ 
is in $E(\Q^\ab)$. 

Moreover, there exist finitely many primes, say $p_1, \ldots, p_k$, such that $\rho_1, \ldots,  \rho_m, \eta_1, \ldots, \eta_n$ 
are contained in the field $F = \Q(\mu_{p_1^{\infty}}, \ldots, \mu_{p_k^{\infty}})$. 
Let $\gamma_{p_1^{s}}$ be a primitive $p_1^s$-th root of unity for any integer $s \ge 1$. 
Then, for any integer $s \ge 1$, the point 
$$
(\rho_1\gamma_{p_1^{s}}^{k_1}+ \cdots + \rho_m \gamma_{p_1^{s}}^{k_m}, \eta_1\gamma_{p_1^{s}}^{l_1} + \cdots + \eta_n \gamma_{p_1^{s}}^{l_n})
$$ 
is in fact in $E(F)$. 
Clearly, these points are of uniformly bounded height and form an infinite set. 
However, from Example~\ref{exam:Zp} $E(F)$ has Property $\N$. 
We then obtain a contradiction. 
The desired result now follows.

\subsection{Proof of Theorem~\ref{thm:Ptor}} 
By contradiction, assume that $E(F)_\tor$ is an infinite set. 
Then, since $E$ is defined over $F$, $E(F)_\tor$ is in fact an infinite abelian group.
We claim that there exist an infinite subgroup $G$ of $E(F)_\tor$ and a prime $p$ such that $p$ does not divide the order of any element of $G$. 

Indeed, if the orders of all the torsion points in $E(F)_\tor$ have prime factors that belong to a finite set of primes, 
then we  choose a prime $p$ which is not a prime of such set and $G = E(F)_\tor$.  
If otherwise the orders of all the torsion points in $E(F)_\tor$ have infinitely many distinct prime factors, say $p_1, p_2, \ldots$, 
then we choose $p=p_1$ and the subgroup $G$ generated by all the $p_i$-torsion points in $E(F)_\tor$, $i = 2, 3, \ldots$.

Now, due to the choices of $G$ and $p$, we have $pG= G$. 
This contradicts with the assumption that $E(F)$ has the property (P). 
So, $E(F)_\tor$ must be finite.
This completes the proof.

\subsection{Proof of Theorem~\ref{thm:preper}}
We essentially use the strategy in the proof of \cite[Proposition 3.1]{DZ}. 
It suffices to consider the case when $\cS$ is an infinite set. 

Let $\Omega$ be the set of preperiodic points in $\cS$ of the endomorphism $f$. 
Let $\Omega_0$ be the set of periodic points in $\Omega$. 
By contradiction, assume that $\Omega$ is infinite. 

If $\Omega_0$ is infinite, then since $f(\Omega_0) = \Omega_0$, we obtain a contradiction with the fact that $\cS$ has the property (P). 
So, $\Omega_0$ must be finite. 

For any point $P \in \cS$, we define the set 
$$
O_f^{-}(P) = \{Q \in \cS \setminus \Omega_0:\, Q \ne P, f^{(n)}(Q) = P \textrm{ for some $n > 0$}\}.
$$
Since $\Omega$ is infinite and $\Omega_0$ is finite, we must have that 
there exists a periodic point $P_0 \in \Omega_0$ such that $O_f^{-}(P_0)$ is infinite. 
Note that under the endomorphism $f$, $P_0$ has at most $\deg f$ preimages. 
So, there exists a point $P_1 \in O_f^{-}(P_0)$ such that $f(P_1) = P_0$ and $O_f^{-}(P_1)$ is infinite. 
Repeating this process, we obtain an infinite sequence of points: $P_0, P_1, P_2, ...$ such that $f(P_i) = P_{i-1}$ for any $i >0$. 
It is easy to see that in this sequence all the points are pairwise distinct 
(because the points $P_1, P_2, ...$ are not periodic points). 
Hence, the set $\Gamma = \{P_0, P_1, P_2, ...\} \cup \{P_0, f(P_0), f^{(2)}(P_0), ...\} $ is infinite, 
and $\Gamma \subseteq \cS$ (due to $f(\cS) \subseteq \cS$). 
Clearly, $f(\Gamma) = \Gamma$, 
which violates the property (P) of $\cS$. This completes the proof.

\subsection{Proof of Theorem~\ref{thm:NimpliesP}}
We use the strategy in the proof of \cite[Theorem 3.1]{DZ}. 
It suffices to consider the case when $\cS$ is an infinite set. 

Let $f$ be an endomorphism of $E$ defined over $\Q(\cS)$ with degree $d >1$, 
and let $\Gamma \subseteq \cS$ be a subset such that $f(\Gamma) = \Gamma$. 
Our task is to prove that $\Gamma$ is a finite set. 

By \eqref{eq:f-h}, there exists a number $B = B(E,f) > 0$ such that $\h(f(P)) \ge d\h(P)-B$ for every point $P \in E(\oQ)$. 
We deduce that if $\h(P) \geq 2B$,  then $\h(f(P)) \ge g \h(P) \geq 2B,$ where $g= d -(1/2) \geq 3/2$.
 By iteration we see that if $\h(P) \geq 2B$, then $\h(f^{(n)}(P)) \ge g^n \h(P)$ for any integer $n > 0$.

Now, let $Q \in \Gamma$. Since $f(\Gamma) = \Gamma$, we can construct an infinite sequence $Q_0 = Q, Q_1, Q_2, ...,$ such that $f(Q_i)=Q_{i-1}$ for any $i > 0$. 
We claim that $\h(Q_n) < 2B$ for all large enough $n$. In fact, if $\h(Q_n) \geq 2B$, 
we deduce from the above (taking $P = Q_n$) that $\h(Q) \ge g^n \h(Q_n) \geq 2g^n B$, which cannot hold for large $n$. 
Hence, since $\cS$ has the Northcott property, there are only finitely many distinct points in the above sequence $Q_0, Q_1, Q_2, ...$. 
So, there exist $r < s$ arbitrarily large and such that $Q_r = Q_s$, 
and in fact this sequence forms a cycle of $f$ (due to the construction of this sequence).   
This means that $Q = Q_0$ belongs to the period $\{Q_s, Q_{s-1},...,Q_r \}$ and in particular $\h(Q)< 2B$.
Hence, $\Gamma$ consists of elements in $\cS$ with height at most $2B$, 
and is therefore a finite set (due to the Northcott property of $\cS$), as asserted.

\subsection{Proof of Theorem~\ref{thm:critP}}
It suffices to consider the case when $\cS$ is an infinite set. 

By contradiction, suppose that $\cS$ does not have the property (P). 
Then, by definition, there is an infinite subset $\Gamma$ of $\cS$ such that 
$f(\Gamma) = \Gamma$ for some endomorphism $f$ of degree $d > 1$. 
Since torsion points are mapped to torsion points under $f$ and $\cS$ has only finitely many torsion points, 
without loss of generality we assume that $\Gamma$ only contains non-torsion points. 

Note that $\Gamma$ also does not have the property (P). 
Then, by Theorem~\ref{thm:NimpliesP}, $\Gamma$ does not have the Northcott property. 
So, we can choose an infinite sequence of 
pairwise distinct non-torsion points in $\Gamma$, say $P_1, P_2, ...$, such that 
$$
\hh(P_i) \le B
$$ 
for each $i \ge 1$ and some constant $B > 0$. 
 
Now, since $f(\Gamma) = \Gamma$, we have that for each integer $i \ge 1$, $\Gamma$ contains a subset of the form 
$$
S_i = \{P_{i1}=P_i, P_{i2}, P_{i3}, ...\}  \subseteq  \Gamma, 
$$
where $f(P_{ij}) = P_{i,j-1}$ for any $j > 1$. 
By \eqref{eq:f-hh}, we obtain that for each $i,j\ge 1$
\begin{equation}  \label{eq:hhP1j}
\hh(P_{ij}) = \frac{1}{d^{j-1}} \hh(P_{i}) \le \frac{B}{d^{j-1}}.
\end{equation}

Suppose that the set $S_1$ is an infinite set. 
Then, by \eqref{eq:hhP1j} and noticing $d > 1$, we have that $\hh(P_{1j}) \to 0$ when $j \to \infty$. 
However, by assumption, the canonical heights of the non-torsion points in $\cS$ have a uniform positive lower bound. 
Hence, $S_1$ must be a finite set. 
Similarly, we can see that all such subsets $S_i$ are finite sets. 
Hence, by construction each $S_i$ is in fact a cycle of $f$. 

If $\max_{i \ge 1} |S_i|$ is not finite, without loss of generality we assume $|S_1| < |S_2| < ...$. 
Write $m_i = |S_i|$ for each $i \ge 1$. Then, $m_i \to \infty$ when $i \to \infty$. 
By \eqref{eq:hhP1j}, we have for each $i \ge 1$, 
$$
\hh(P_{im_i}) \le \frac{B}{d^{m_i -1}}.
$$
So, $\hh(P_{im_i}) \to 0$ when $i \to \infty$. 
We again obtain a contradiction with the positive lower bound assumption for the canonical heights. 
Hence, $\max_{i \ge 1} |S_i|$ must be  finite.
     
Note that each $S_i$ is a cycle of $f$. 
Then, since $\max_{i \ge 1} |S_i|$ is  finite,
 there exists a positive integer $n$ such that $f^{(n)}(P_i) = P_i$ for each $i \ge 1$. 
That is, the endomorphism $f^{(n)} - 1$ has an infinite kernel. 
So, it is the zero morphism, and thus $f^{(n)} = 1$. 
This means that $f$ is an automorphism and so $\deg f = 1$, 
which contradicts with $\deg f > 1$. 
We then complete the proof.

\section*{Acknowledgement}
The authors would like to thank the referee for his/her careful reading and  valuable comments. 
The authors are grateful to Igor Shparlinski and Joseph Silverman for helpful discussions, 
and thank Lukas Pottmeyer for valuable comments (especially for Example~\ref{exam:pres}).  
For the research, the first author was supported by the Australian Research Council Grant DP180100201, by NSERC, and by the Max Planck Institute for Mathematics.
The second author was supported by the Australian Research Council Grant DE190100888 
and by the Natural Science Foundation of Guangdong Province.

\end{document}